\documentclass[oneside,12pt]{amsart}          % please do not change
\usepackage[a4paper]{geometry}	    
\usepackage{amsfonts,amsmath,latexsym,amssymb,enumerate,hyperref,ifthen} % 
\usepackage[mathscr]{eucal}
\usepackage[utf8]{inputenc}
\usepackage[T1]{fontenc}

\def\N{\mathbb{N}}
\def\R{\mathbb{R}}
\def\F{\mathbb{F}}
\def\Q{\mathbb{Q}}

\def\intr{\mathop{\mbox{\rm int}}}
\def\eps{\varepsilon}
\def\phi{\varphi}

\def\twolinelim#1{\lim_{\hbox{\tiny$\begin{array}{c}#1\end{array}$}}}

\newtheorem{theorem}{Theorem}
\def\Thm#1#2{\ifthenelse{\equal{#1}{*}}{\begin{theorem*}#2\end{theorem*}}
  {\begin{theorem}\label{T#1}#2\end{theorem}}}
\newtheorem{Atheorem}{Theorem}

\def\thm#1{Theorem~\ref{T#1}}

\newtheorem{proposition}[theorem]{Proposition}
\def\Prp#1#2{\ifthenelse{\equal{#1}{*}}{\begin{proposition*}#2\end{proposition*}}
             {\begin{proposition}\label{P#1}#2\end{proposition}}}
\def\prp#1{Proposition~\ref{P#1}}
\newtheorem{corollary}[theorem]{Corollary}
\def\Cor#1#2{\ifthenelse{\equal{#1}{*}}{\begin{corollary*}#2\end{corollary*}}
             {\begin{corollary}\label{C#1}#2\end{corollary}}}

\newtheorem{lemma}[theorem]{Lemma}
\def\Lem#1#2{\ifthenelse{\equal{#1}{*}}{\begin{lemma*}#2\end{lemma*}}
             {\begin{lemma}\label{L#1}#2\end{lemma}}}

\newtheorem{example}[theorem]{Example}
\long\def\Exa#1#2{\ifthenelse{\equal{#1}{*}}{\begin{example*}\rm #2\end{example*}}
            {\begin{example}\label{Ex#1}\rm #2\end{example}}}

\theoremstyle{definition}
\newtheorem{definition}[theorem]{Definition}
\def\Defi#1#2{\ifthenelse{\equal{#1}{*}}{\begin{definition*}#2\end{definition*}}
      {\begin{definition}\label{D#1}#2\end{definition}}}

\newtheorem{remark}[theorem]{Remark}
\def\Rem#1#2{\ifthenelse{\equal{#1}{*}}{\begin{remark*}#2\end{remark*}}
             {\begin{remark}\label{R#1}#2\end{remark}}}

\def\eq#1{{\rm(\ref{E#1})}}
\def\Eq#1#2{\ifthenelse{\equal{#1}{*}}
  {\begin{equation*}\begin{aligned}#2\end{aligned}\end{equation*}}
  {\begin{equation}\begin{aligned}\label{E#1}#2\end{aligned}\end{equation}}}

\begin{document}
\begin{flushright}
Math. Inequal. Appl. \textbf{15}(2) (2012), 289--299.\\
\href{http://dx.doi.org/10.7153/mia-15-24}{doi: 10.7153/mia-15-24} \\[1cm]
\end{flushright}

\title[]{On strong $(\alpha,\F)$-convexity}
\author{Judit Mak\'o, Kazimierz Nikodem, and Zsolt P\'ales}
\address{Judit Mak\'o, Institute of Mathematics, University of Debrecen,
H-4010 Debrecen, Pf.\ 12, Hungary}
\email{makoj@science.unideb.hu}

\address{Kazimierz Nikodem, Department of Mathematics and Computer Science, University of
Bielsko-Biała, ul.\ Wil\-lo\-wa~2, 43-309 Bielsko-Biała, Poland}
\email{knikodem@ath.bielsko.pl}

\address{Zsolt P\'ales, Institute of Mathematics, University of Debrecen,
H-4010 Debrecen, Pf.\ 12, Hungary}
\email{pales@science.unideb.hu}


\thanks{This research has been supported by the Hungarian
Scientific Research Fund (OTKA) Grant NK81402 and by the
TÁMOP 4.2.1./B-09/1/KONV-2010-0007 project implemented through the New
Hungary Development Plan co-financed by the European Social Fund, and the
European Regional Development Fund.}

%\dedicatory{}

\begin{abstract}
In this paper, strongly $(\alpha,T)$-convex functions, i.e., functions
$f:D\to \R$ satisfying the functional inequality
\Eq{*}
{
f(tx+(1-t)y)\leq tf(x)+(1-t)f(y)-t\alpha\big((1-t)(x-y)\big)-(1-t)\alpha\big(t(y-x)\big)
}
for $x,y\in D$ and $t\in T\cap[0,1]$ are investigated. Here $D$ is a convex set in a linear 
space, $\alpha$ is a nonnegative function on $D-D$, and $T\subseteq\R$ is a nonempty set.
The main results provide various characterizations of strong $(\alpha,T)$-convexity 
in the case when $T$ is a subfield of $\R$.
\end{abstract}

\maketitle

\section{Introduction}

Let $I\subset \R$ be an interval and $c$ be a positive number. A
function $f:I\to \R$ is called \emph{strongly convex with modulus
$c$} if
\Eq{*}{
 f(tx+(1-t)y)\leq t f(x)+(1-t)f(y)-ct(1-t)(x-y)^2,
}
for all $x,y\in I$ and $t\in [0,1]$. Strongly convex functions
have been introduced by Polyak \cite{Pol66} and they play an
important role in optimization theory and mathematical economics.
Many properties of them can be found in the literature (see, e.g.
\cite{Pol96}, \cite{RobVar73}, \cite{Via82}, \cite{MerNik11}). It is known,
for instance, that a function $f:I\to \R$ is strongly convex with
modulus $c$ if and only if for every $x_0\in\intr I$ there exists an
$a\in \R$ such that
\Eq{*}{
f(x)\geq c(x-x_0)^2 + a(x-x_0) + f(x_0), \qquad x\in I,
}
i.e., $f$ has a quadratic support at $x_0$. If $f$ is
differentiable and strongly convex with modulus $c$ then its
derivative $f'$ is "strongly monotone" in the sense:
$(f'(x)-f'(y))(x-y)\geq 2c(x-y)^2$, $x,y\in I$ \
(cf.\ \cite[p. 268]{RobVar73}). \vspace{1ex}

In this paper we introduce the class of strongly $(\alpha,T)$-convex
functions (which is much wider then the class of strongly convex
functions) and present, among other, some generalizations of the
results mentioned above.
\vspace{1ex}

Let $X$ be a real linear space.
For a nonempty convex subset $D\subseteq X$, denote
\Eq{*}
{
  D^*:=D-D:=\{x-y:x,y\in D\}.
}
Given a nonnegative \textit{even} function $\alpha:D^*\to \R_+$ and a nonempty set $T\subseteq\R$
such that $T\cap[0,1]$ is nonempty,
we say that a map $f:D\to \R$ is \textit{strongly $(\alpha,T)$-convex}, if for all $x,y\in D$
and $t\in T\cap[0,1]$,
\Eq{conv}
{
f(tx+(1-t)y)\leq tf(x)+(1-t)f(y)-t\alpha\big((1-t)(x-y)\big)-(1-t)\alpha\big(t(y-x)\big)
}
holds.
If \eq{conv} holds with $T=\{1/2\},$ i.e., for all
$x,y\in D,$
\Eq{Jconv}
{
f\big(\tfrac{x+y}{2}\big)\leq \tfrac{f(x)+f(y)}{2}-\alpha\big(\tfrac{x-y}{2}\big),
}
then the function $f$ is called \textit{strongly $\alpha$-Jensen convex.}
If $T\supseteq[0,1]$, then the function $f$ is termed \textit{strongly $\alpha$-convex.}
By the nonnegativity of $\alpha$, we have that strongly
$\alpha$-Jensen convex, strongly $(\alpha,T)$-convex, and strongly $\alpha$-convex
functions are always convex in the same sense, respectively. More generally,
if $\alpha,\beta:D^*\to\R_+$ and $\alpha\geq\beta$, then strong convexity with respect to
$\alpha$ in some sense implies strong convexity with respect to $\beta$ in the same sense.
Note also that for $\alpha (x) = cx^2$ strong $\alpha$-convexity coincides with strong
convexity with modulus $c$.

\section{Strengthening the strong Jensen convexity}

In the next theorem and corollary, which are particular cases of the theorem
in \cite{MakPal11a}, the strong $\alpha$-Jensen convexity property will be strengthened.
We provide their proof here because it is much simpler and more transparent than in the
general case.

\Thm{S}{Let $f:D\to\R$ be a strongly $\alpha$-Jensen convex function.
Then $f$ is strongly $\widetilde{\alpha},$-Jensen convex on $D$, where
\Eq{S}
{
\widetilde{\alpha}(u):=\sup\big\{n^2\alpha\big(\tfrac{u}{n}\big)\mid n\in\N\big\}
  \qquad(u\in D^*).
}}

\begin{proof}
Assume that $f:D\to\R$ is strongly $\alpha$-Jensen convex and let $n\in\N.$ Let $x,y\in D.$
Consider the segment $[x,y]$ and divide it into $2n$ pieces of equal subsegments. For this, we
define the system of points $x_0=x$, $x_1,\dots,x_{2n-1}$, $x_{2n}=y$ in the following way
\Eq{xi}{
  x_{i}:=x+\tfrac{i}{2n}(y-x)
  \qquad(i\in\{0,\dots,2n\}).
}
We have the following two obvious identities:
\Eq{*}{
  x_{i+1}-x_{i-1}&=\tfrac{y-x}{n} \qquad\mbox{and}\qquad
  x_{i}=\tfrac{x_{i-1}+x_{i+1}}{2}.
    \qquad(i\in\{1,\dots,2n-1\}).
}
Therefore, by the strong $\alpha$-Jensen convexity of $f$, we get
\Eq{c}{
  f(x_{i})&\leq \tfrac{f(x_{i-1})
    +f(x_{i+1})}{2}-\alpha\big(\tfrac{x-y}{2n}\big)
}
for all $i\in\{1,\dots,2n-1\}.$
Multiplying \eq{c} by
\Eq{*}
{
i,&\quad \mbox{if}\quad &&i\in \{1,\dots,n\}\\
2n-i,&\quad \mbox{if}\quad &&i\in \{n+1,\dots,2n-1\}
}
and adding the inequalities so obtained, we get that
\Eq{EEE}{
  &\sum_{i=1}^{n}if(x_{i})
  +\sum_{i=n+1}^{2n-1}(2n-i)f(x_i)\\
  &\leq\sum_{i=1}^n \tfrac{i}2(f(x_{i-1})+f(x_{i+1}))
  +\sum_{i=n+1}^{2n-1}\tfrac{2n-i}2(f(x_{i-1})+f(x_{i+1}))
    -\bigg(\sum_{i=1}^n i+\sum_{i=n+1}^{2n-1}(2n-i)\bigg)\alpha\big(\tfrac{x-y}{2n}\big).
}
The coefficient of $f(x_i)$ in inequality \eq{EEE} is the following:
\Eq{*}
{
i-\tfrac{i+1}2-\tfrac{i-1}2&=0,\qquad \mbox{if} \qquad 1\leq i\leq n-1,\\
n-\tfrac{n-1}2-\tfrac{2n-(n+1)}2&=1, \qquad \mbox{if} \qquad i=n,\\
(2n-i)-\tfrac{2n-i-1}2-\tfrac{2n-i+1}2&=0,\qquad \mbox{if} \qquad n+1\leq i\leq 2n-1.
}
The terms $f(x_0)=f(x)$ and $f(x_{2n})=f(y)$ appear only on the right hand side
of \eq{EEE} with coefficients
\Eq{*}{
  \tfrac{1}{2}
  \qquad\mbox{and} \qquad
  \tfrac{2n-(2n-1)}{2}=\tfrac{1}{2},
}
respectively. Finally, the coefficient of the error function is
\Eq{*}{
\sum_{i=1}^{n} i+\sum_{i=n+1}^{2n-1}(2n-i)=\tfrac{n(n+1)}{2}+\tfrac{(n-1)n}2=n^2.
}
Thus, \eq{EEE} reduces to,
\Eq{*}
{
f\big(\tfrac{x+y}{2}\big)=f(x_n)
   \leq \tfrac{f(x_0)+f(x_{2n})}2-n^2\alpha\big(\tfrac{x-y}{2n}\big)
   = \tfrac{f(x)+f(y)}2-n^2\alpha\big(\tfrac{x-y}{2n}\big).
}
Therefore, we get that $f$ is strongly $\widetilde{\alpha}$-Jensen convex, where
$\widetilde{\alpha}$ is defined by \eq{S}. This completes the proof.
\end{proof}

\Exa{1}{Assume that a function $f:I\to\R$ is strongly $\sin^2$-Jensen convex, i.e.,
\Eq{*}
{
f\big(\tfrac{x+y}{2}\big)\leq \tfrac{f(x)+f(y)}{2}-\sin^2\big(\tfrac{x-y}{2}\big)\qquad(x,y\in I),
}
then, by \thm{S}, $f$ is also $\widetilde{\sin^2}$-Jensen convex, where, for $u\in\R$,
\Eq{*}{
\widetilde{\sin^2}(u)=\sup\{n^2\sin^2(\tfrac{u}{n})\mid n\in\N\}
= u^2\sup\{\tfrac{n^2}{u^2}\sin^2(\tfrac{u}{n})\mid n\in\N\}
= u^2 \lim_{n\to\infty}\tfrac{\sin^2\big(\tfrac{u}{n}\big)}{\big(\tfrac{u}{n}\big)^2}=u^2.
}
This means that $f$ is also strongly Jensen convex with modulus $1$, i.e.,
\Eq{*}
{
f\big(\tfrac{x+y}{2}\big)\leq \tfrac{f(x)+f(y)}{2}-\big(\tfrac{x-y}{2}\big)^2\qquad(x,y\in I).
}}

\Cor{S}{Assume that there exists
$u\in \frac12D^*\setminus\{0\}$ such that
\Eq{lim}{
\limsup_{n\to\infty}n^2\alpha\big(\tfrac{u}{n}\big)=\infty.
}
Then there is no strongly $\alpha$-Jensen convex function on $D$.
}

\begin{proof}
Assume that there exists a strongly $\alpha$-Jensen convex function $f:D\to\R$ and 
let $u\in \frac12D^*\setminus\{0\}$ such that \eq{lim} holds. By \thm{S},
we get that $f$ is $\widetilde{\alpha}$-Jensen convex on $D$.
The definition of $\widetilde{\alpha}$ implies that,
\Eq{leq,}{
\widetilde{\alpha}(u)\geq\limsup_{n\to\infty}n^2 \alpha\big(\tfrac{u}{n}\big)=\infty.
}
Hence, we get that $-\widetilde{\alpha}(u)=-\infty,$ which means that the error term in 
\eq{Jconv} is equal to $-\infty$ for some $x,y\in D$. This is a contradiction resulting the statement.
\end{proof}

\Rem{1}{If $X$ is a normed space, $\eps$, $p$ are positive constants, and
$\alpha(u):=\eps\|u\|^p$ for $u\in D^*$, then condition \eq{lim} holds if and only if
$p<2$.}

\section{On strong $(\alpha,\F)$-convexity}

Let $\F$ be a subfield of $\R$.
Given a real linear space $X$, a function $\phi:X\to\R$ is called \textit{$\F$-linear},
if it is \textit{additive}, i.e., for all $x,y\in X,$
\Eq{*}{
\phi(x+y)=\phi(x)+\phi(y),
}
and it is \textit{$\F$-homogeneous}, i.e., for all $x\in X$ and for all $\lambda\in\F,$
\Eq{*}{
\phi(\lambda x)=\lambda\phi(x).
}
As it is well-known, additive functions on $X$ are automatically $\Q$-linear.

A function $\psi:X\to\R$ is called \textit{$\F$-sublinear}, if it is \textit{subadditive},
i.e., for all $x,y\in X$,
\Eq{*}{
\psi(x+y)\leq\psi(x)+\psi(y),
}
and it is \textit{positively $\F$-homogeneous}, i.e., for all $x\in X$ and for all
$\lambda\in\F\cap\R_+,$
\Eq{*}{
\psi(\lambda x)=\lambda\psi(x).
}
The \textit{$\F$(-algebraic) dual} $X'_\F$ of the space $X$ is defined by
\Eq{*}{
  X'_\F:=\{\phi:X\to\R\mid \phi \mbox{ is $\F$-linear}\}.
}
As it is well-known, additive functions are automatically $\Q$-linear, therefore
we always have the inclusion $X'_\F\subseteq X'_\Q$.

The following result is a consequence of the standard separation/sandwich theorems
(cf.\ Mazur--Orlicz \cite{MazOrl53}, Holmes \cite{Hol75}, Nikodem--Páles--Wąsowicz
\cite{NikPalWas99}) or of the Rodé theorem (\cite{Rod78}).

\Thm{HB}{Let $\psi:X\to \R$ be an $\F$-sublinear function, then there exists a
$\phi\in  X'_\F,$ such that $\phi\leq\psi$.}

Let $D\subseteq X$ be a nonempty set. We say that $D$ is \textit{$\F$-algebraically open}
if, for all $x_0\in D$ and $h\in X,$ there exits an $\eps\in]0,+\infty[$ such that
$[x_0,x_0+th]\subset D$ for all $t\in [0,\eps]\cap \F.$
In what follows, assume that $D\subseteq X$ is a nonempty $\F$-algebraically open convex set.

\Prp{H}{[Cf.\ \cite{BorPal06}] Let $f:D\to \R$ be an $\F$-convex function,
$x_0\in D$ and $h\in X,$ then the mapping
\Eq{*}
{
t\mapsto \frac{f(x_0+th)-f(x_0)}{t}
}
is nondecreasing on the set $\{t\in\F\setminus\{0\}\mid x_0+th\in D\}$.
}

The generalized \textit{$\F$-directional derivative} of an $\F$-convex function $f:D\to\R$ at
$x_0\in D$ in a direction $h\in X$, denoted by $f'_{\F}(x_0,h)$, is defined as follows
\Eq{*}{
f'_{\F}(x_0,h):=\twolinelim{t\to 0^+ \\ t\in\F}\frac{f(x_0+th)-f(x_0)}{t}.
}
Note that this generalized $\F$-directional derivative has similar properties as
the standard directional derivative.

\Prp{D1}{[Cf.\ \cite{BorPal06}] Let $f:D\to\R$ be an $\F$-convex function and $x_0\in D$ be an
arbitrary element of $D$, then the mapping $h\mapsto f'_{\F}(x_0,h)$ is $\F$-sublinear.}

Our main result is contained in the following theorem.

\Thm{E}{For any function $f:D\to\R$, the following conditions are equivalent:
\begin{enumerate}[(i)]
  \item $f$ is strongly $(\alpha,\F)$-convex.
  \item $f$ is $\F$-directionally differentiable at every point of $D$,
     and for all $x_0\in D,$ the map $h\mapsto f'_\F(x_0,h)$ is $\F$-sublinear on $X$,
     furthermore for all $x_0,x\in D,$
    \Eq{f'}{
     f(x)\geq f(x_0)+f'_{\F}(x_0,x-x_0)+\alpha(x-x_0).
     }
  \item For all $x_0\in D$, there exits an element $\phi\in X'_{\F}$ such that
    \Eq{phi}{
     f(x)\geq f(x_0)+\phi(x-x_0)+\alpha(x-x_0)
     \quad \mbox{for all} \quad x\in D.
     }
\end{enumerate}}

\begin{proof}
$(\mbox{i}) \Rightarrow (\mbox{ii})$ It is evident that $f$ is
$\F$-convex, which implies that $f$ is $\F$-directionally differentiable at
every point of $D$, moreover for all $x_0\in D$ the map $h\mapsto
f'_\F(x_0,h)$ is $\F$-sublinear on $X$. To prove \eq{f'}, let $x_0,x\in D$
be arbitrary. Since $D$ is $\F$-algebraically open, there exists
an $\eps\in ]0,1[$ such that $x+\frac{t}{1-t}(x-x_0)=x_0+\frac{x-x_0}{1-t}\in D$
for all $t\in[0,\epsilon]\cap\F$. Let $h:=\frac{x-x_0}{1-t}.$
By the strong $(\F,\alpha)$-convexity of $f$, we get that
\Eq{*} {
f(x_0+th) \leq
(1-t)f(x_0)+tf(x_0+h)-t\alpha\big((1-t)h\big)-(1-t)\alpha(-th)
  \qquad \mbox{for all}\quad t\in [0,\eps]\cap\F.
}
Rearranging the above inequality, we have that
\Eq{*} {
  f(x_0+h)\geq
               f(x_0)+\frac{f(x_0+th)-f(x_0)}{t}
               +\alpha\big((1-t)h\big)+\frac{1-t}{t}\alpha(-th)
        \quad \mbox{for all}\qquad t\in ]0,\eps]\cap\F.
}
By the nonnegativity of $\alpha,$
\Eq{*} {
 f(x_0+h) \geq f(x_0)+\frac{f(x_0+th)-f(x_0)}{t}+\alpha\big((1-t)h\big)\qquad \mbox{for all}\quad t\in ]0,\eps]\cap\F.
}
Substituting $h=\frac{x-x_0}{1-t},$ we obtain that, for all $t\in]0,\epsilon]\cap\F$,
\Eq{e}{
f\Big(x_0+\frac{x-x_0}{1-t}\Big)\geq
f(x_0)+\frac{1}{1-t}\frac{f\big(x_0+\frac{t}{1-t}(x-x_0)\big)-f(x_0)}{t/(1-t)}
+\alpha(x-x_0)
}
holds.
By the $\F$-convexity of $f$, the mapping $s\mapsto f\big(x_0+s(x-x_0)\big)$ is 
continuous on $[0,\epsilon]\cap\F$, whence we get
\Eq{*}{
\twolinelim{t\to 0^+\\t\in\F}f\Big(x_0+\frac{x-x_0}{1-t}\Big)=f(x),
}
furthermore, the limit
\Eq{*}{
\twolinelim{t\to 0^+\\t\in\F}\frac{f\big(x_0+\frac{t}{1-t}(x-x_0)\big)-f(x_0)}{t/(1-t)}
}
exists and equals $f'_\F(x_0,x-x_0).$ Thus, taking the limit $t\to 0^+$ for $t\in\F$ in \eq{e},
we get \eq{f'}, which completes the proof of (ii).

$(\mbox{ii}) \Rightarrow (\mbox{iii})$
Assume that $f$ is $\F$-directionally differentiable at every point of $D$ and for 
all $x_0\in D,$ $h\mapsto f'_{\F}(x_0,h)$
is $\F$-sublinear. By \thm{HB}, there exists an element
$\phi\in X'_{\F},$ such that
\Eq{*}
{
f'_{\F}(x_0,h)\geq \phi(h) \qquad \mbox{for all} \quad h\in X.
}
This and \eq{f'} implies that \eq{phi} holds.

$(\mbox{iii}) \Rightarrow (\mbox{i})$ Let $x,y\in D,$ $t\in [0,1]\cap\F$, and set $x_0:=tx+(1-t)y$.
Then, by $(\mbox{iii})$, we have
\Eq{*}
{
f(x)&\geq f(tx+(1-t)y)+\phi\big((1-t)(x-y)\big)+\alpha\big((1-t)(x-y)\big),\\
f(y)&\geq f(tx+(1-t)y)+\phi\big(t(y-x)\big)+\alpha\big(t(y-x)\big).
}
Multiplying the first inequality by $t$ and the second inequality by $1-t$
and adding up the inequalities so obtained, we get $\eq{conv}.$
\end{proof}

\Cor{alpha}{Let $f:D\to \R$ be a strongly $(\alpha,\F)$-convex function.
Then $\alpha(0)=0$, $\alpha$ is $\F$-directionally differentiable at $0$ and
\Eq{*}{
\alpha'_{\F}(0,h)=0\qquad (h\in X).}
}

\begin{proof} If $f$ is strongly $(\alpha,\F)$-convex then property (ii) of \thm{E} holds.
Let $x_0\in D$ be fixed and $h\in X.$ Since $D$ is $\F$-algebraically open, there exists
an $\eps\in ]0,1[\cap\F$ such that $x_0+\eps h\in D$. Then, substituting $x=x_0+th$ into
\eq{f'} (where $t\in]0,\eps[\cap\F$), we get
\Eq{f*}{
     f(x_0+th)\geq f(x_0)+f'_\F(x_0,th)+\alpha(th).
     }
Taking $h=0$, by the nonnegativity of $\alpha$, it follows that $\alpha(0)=0$.
On the other hand, rearranging the above inequality,
\Eq{*}{
  \frac{f(x_0+th)-f(x_0)}{t}\geq f'_\F(x_0,h)+\frac{\alpha(th)}{t}.
}
By taking the limit $t\to0+$ and using the nonnegativity of $\alpha$ again, it follows that
\Eq{++}{
  \twolinelim{t\to 0^+\\t\in\F}\frac{\alpha(th)}{t}=0.
}
Therefore, $\alpha'_\F(0,h)=0$, proving that the $\F$-directional derivative of 
$\alpha$ at the origin is zero.
\end{proof}

In the following result, we strengthen $(\alpha,\F)$-convexity.

\Cor{S3}{Let $f:D\to \R$ be a strongly $(\alpha,\F)$-convex function. Then
$f$ is strongly $(\widehat{\alpha},\F)$-convex, where
$\widehat{\alpha}:D^*\to \R$ is defined by
\Eq{S3a}{
  \widehat{\alpha}(u):=\sup\big\{\tfrac{\alpha(tu)}{t}\mid t\in ]0,1]\cap\F\big\}.
}}

\begin{proof}
Let $f:D\to\R$ be a strongly $(\alpha,\F)$-convex function and
$x_0,x$ be arbitrary elements of $D$. By \prp{H}, the mapping $t\mapsto
\frac{f(x_0+t(x-x_0))-f(x_0)}{t}$ is nondecreasing on $]0,1]\cap\F$. Thus, for
all $t\in]0,1]\cap\F,$ we have that
\Eq{*}{
f(x)-f(x_0)\geq \frac{f(x_0+t(x-x_0))-f(x_0)}{t}.
}
By \thm{E}, we get that
\Eq{*}{
f(x_0+t(x-x_0))-f(x_0)\geq f'_\F(x_0,t(x-x_0))+\alpha(t(x-x_0))\qquad (t\in ]0,1]\cap\F).
}
Combining the above inequalities, we obtain
\Eq{*}{
f(x)-f(x_0)\geq f'_\F(x_0,x-x_0)+\tfrac{\alpha(t(x-x_0))}{t}\qquad (t\in ]0,1]\cap\F).
}
Using \thm{E}, this means that the function $f$ is also strongly
$(\widehat{\alpha},\F)$-convex, which completes the proof.
\end{proof}

\section{The $\F$-subgradient of $(\alpha,\F)$-convex functions}

For any $f:D\to \R$ function and $x_0\in D$, define the
\textit{$\F$-subgradient} of $f$ at $x_0$ by
\Eq{p}
{
\partial_{\F}f(x_0):=\{\phi\in X'_\F\mid f(x)\geq f(x_0)+\phi(x-x_0)
\quad \mbox{for all} \quad x\in D\}.
}
Obviously, $\partial_{\F}f(\cdot)$ can be considered as a set-valued mapping defined on $D$
with values in $2^{X'_\F}$.

For $\F$-convex functions, the $\F$-subdifferential $\partial_{\F}f(x_0)$ can also
be expressed in terms of the $\F$-directional derivative of $f$ at $x_0$.

\Prp{3}{[Cf.\ \cite{BorPal06}] Let $f:D\to\R$ be an $\F$-convex function. Then, for all $x_0\in D$,
\Eq{ps}
{
\partial_{\F}f(x_0)=\{\phi\in X'_\F\mid f'_\F(x_0,h)\geq \phi(h)\quad \mbox{for all} \quad h\in X\}.
}}

To describe the properties of the $\F$-subdifferential of
strongly $(\alpha,\F)$-convex functions, we need to recall and define certain
generalized monotonicity concepts whose original versions were introduced 
by Minty \cite{Min64b} and R.~T.~Rockafellar (\cite{Roc66c}, \cite{Roc70b}, \cite{Roc70})
in order to characterize the subdifferentials of convex functions.

We say that a set-valued mapping $ \Phi:D \to 2^{X'_\Q}$ is {\em $\alpha$-monotone} if
\Eq{monAB}{
\varphi(y-x) + \alpha(y-x) + \psi(x-y) + \alpha(x-y) \leq 0
}
holds for every $x,y\in D$, $\varphi\in\Phi(x)$, and $\psi\in\Phi(y)$.
We call a set-valued mapping $\Phi:D\to 2^{X'_\Q}$
{\em $\alpha$-cyclically monotone} if the inequality
\Eq{cycleA}{
\sum_{j=0}^{n} \Big(\varphi_j(x_{j+1}-x_j) +\alpha(x_{j+1}-x_j)\Big)\leq 0
}
is fulfilled for every $n\in\N$, $x_j\in D $ $(j\in\{0,1,\dots,n,n+1\})$
with $x_{n+1} = x_0$, and $ \varphi_j\in \Phi (x_j) $ $(j\in\{0,1,\dots,n\})$.
Obviously, by taking $n=2$ in the above definition, $\alpha$-cyclical monotonicity
implies $\alpha$-monotonicity, however, the reversed implication may not be valid.
In the particular case when $\alpha$ is identically zero, we simply speak about
{\em monotone} and {\em cyclically monotone} set-valued maps. Observe that, by the
nonnegativity of $\alpha$, the properties $\alpha$-monotonicity and $\alpha$-cyclical
monotonicity imply monotonicity and cyclical monotonicity, respectively.

We say that a mapping $ \Phi:D \to 2^{X'_\F}$ is {\em $\F$-maximal monotone} if $\Phi$ is
monotone and, for any monotone mapping $ \Psi:D\to 2^{X'_\F}$ fulfilling
$\Phi(x)\subseteq \Psi(x)$ for all $x\in D$, also $\Phi(x)=\Psi(x)$ holds for every $x\in D$.
In particular, $\Q$-maximal monotone mappings are called {\em maximal monotone}.

In the next result we summarize the properties of the $\F$-subdifferential of
strongly $(\alpha,\F)$-convex functions.

\Thm{P}{Let $f:D\to \R$ be an $(\alpha,\F)$-convex function.
Then, for every $x_0\in D$, $\partial_{\F}f(x_0)$ is a nonempty convex subset
in $X'_\F$ which is closed with respect to the pointwise convergence and
\Eq{pa}
{
\partial_{\F}f(x_0)=\{\phi\in X'_\F\mid f(x)\geq f(x_0)+\phi(x-x_0)+\alpha(x-x_0)
\quad \mbox{for all} \quad x\in D\}.
}
Furthermore, the map $\partial_{\F}f:D\to 2^{X'_\F}$ is $\F$-maximal monotone
and $\alpha$-cyclically monotone.}

\begin{proof} The convexity and closedness (with respect to the pointwise convergence)
of $\partial_{\F}f(x_0)$ directly follows from its definition. The inclusion $\supseteq$
in \eq{pa} is a consequence of the definition \eq{p}. To prove the reversed inclusion,
let $\phi\in\partial_{\F}f(x_0)$ be arbitrary. Then, by \prp{3}, $f'_\F(x_0,h)\geq\phi(h)$
holds for all $h\in X$. On the other hand, by the second assertion of \thm{E}, for all $x\in D$,
we have
\Eq{*}{
 f(x)\geq f(x_0)+f'_\F(x_0,x-x_0)+\alpha(x-x_0).
}
Hence, for all $x\in D$,
\Eq{*}{
 f(x)\geq f(x_0)+\phi(x-x_0)+\alpha(x-x_0).
}
This proves that $\phi$ also belongs to the right hand side of \eq{pa} and hence \eq{pa}
holds with equality.

By the third assertion of \thm{E}, the right hand side of \eq{pa} is nonempty, 
which yields the nonemptiness of $\partial_{\F}f(x_0)$.

The $\F$-maximal monotonicity is a consequence of \cite[Theorem 5.4]{BorPal06}.

To prove the $\alpha$-cyclic monotonicity of $f$, let $n\in\N$, $x_j\in D$ for $j\in\{0,1,\ldots,n,n+1\}$
with $x_{n+1} = x_0$, and $\varphi_j\in\partial_{\F}f(x_j)$ for $j\in\{0,1,\ldots,n\}$.
Then, by the third assertion of \thm{E},
\[
    f(x_{j+1})\geq f(x_j)+\phi(x_{j+1}-x_j)+\alpha(x_{j+1}-x_j)
     \qquad (j\in\{0,1,\ldots,n\}).
\]
Adding up these inequalities for $j\in\{0,1,\ldots,n\}$ and using $x_{n+1} = x_0$,
the inequality \eq{cycleA} follows immediately proving the $\alpha$-cyclic monotonicity
of $f$.
\end{proof}

The following statement is analogous to \cite[Theorem 1]{Roc66c}.

\Thm{Tc}{If $\Phi:D\to2^{X'_\F}$ is a nonempty-valued $\alpha$-cyclically monotone
set-valued map, then there exists a strongly $(\alpha,\F)$-convex function $f:D \to\R$
such that $\Phi(x)\subseteq \partial_\F f(x)$ for every $x\in D$.}

\begin{proof}
Let $x_0\in D$ be fixed. For each $x\in D$,
let $S(x)$ denote the set of all finite sums of the form
\[
\sum_{j=0}^{n-1}\Big(\phi_j(x_{j+1}-x_j)+\alpha(x_{j+1}-x_j)\Big),
\]
where $ n\in\N$, $x_j\in D $ $(j\in\{1,\ldots,n\})$ such that
$x_n=x$, and $\phi_j\in\Phi(x_j)$ for $j\in\{0,1,\ldots,n-1\}$.

For any $\phi\in\Phi(x)$, the $\alpha$-cyclical monotonicity of $\Phi$ yields
\[
\sum_{j=0}^{n-1} \Big(\phi_j(x_{j+1}-x_j)+\alpha(x_{j+1}-x_j)\Big)
+\phi(x_0-x)+\alpha(x_0-x)\leq 0.
\]
Thus $ -(\phi(x_0-x)+\alpha(x_0-x))$ is an upper bound for $S(x)$.
Hence, we may define a function $f:D\to\R$ by
\[
f(x) = \sup S(x) \qquad (x\in D).
\]

In order to prove the desired inclusion $\Phi(x)\subseteq \partial_\F f(x)$,
consider arbitrary elements $x,y\in D$, $\phi\in\Phi(x)$, and $\eps > 0$.
Then there exist $n\in\N$, $x_j\in D $ $(j\in\{1,\dots,n\})$ with
$x_n = x$, and $ \phi_j\in\Phi(x_j) $ $(j\in\{0,1,\dots,n-1\})$
such that
\Eq{*}{
\sum_{j=0}^{n-1} \Big(\phi_j(x_{j+1}-x_j)+\alpha(x_{j+1}-x_j)\Big)>f(x)-\eps.
}
The definition of $f(y)$ and the above inequality yields
\Eq{*}{
f(y) &\geq\sum_{j=0}^{n-1} \Big(\phi_j(x_{j+1}-x_j)+\alpha(x_{j+1}-x_j)\Big)
          + \phi(y-x)+\alpha(y-x)\\
     &\geq  f(x)-\eps+\phi(y-x)+\alpha(y-x).
}
Letting $\eps$ tend to $0$, we get
\Eq{+++}{
 f(y) \geq f(x)+\phi(y-x)+\alpha(y-x).
}
Hence, $\phi\in \partial_\F f(x)$ which proves the inclusion
$\Phi(x) \subseteq \partial_\F f(x)$.
In particular, for every $x\in D$, by the nonemptiness of $\Phi(x)$,
there exists an element $\phi\in X'_\F$ such that \eq{+++} holds for
all $y\in D$ showing that the third assertion of \thm{E} is valid.
Thus, by \thm{E}, we obtain that $f$ is strongly $(\alpha,\F)$-convex.
\end{proof}

An immediate consequence of \thm{Tc} is the following result.

\Cor{Ccycmon}{If $\Phi:D\to 2^{X'_\F}$ is
an $\F$-maximal $\alpha$-cyclically monotone mapping,
then there exists a strongly $(\alpha,\F)$-convex function
$ f:D\to\R$ such that $\Phi(x)=\partial_\F f(x)$ for every $x\in D$.}

%\nocite{Kuc85,MroTabTab08,MurTabTab09,TabTab09b,TabTabZol10,BerDoe15,Ger88c}
%\bibliography{publ,funcequ,control}
%\bibliographystyle{amsplain}
%\end{document}

\providecommand{\MR}{\relax\ifhmode\unskip\space\fi MR }

\end{document}